\documentclass[fleqn,10pt,twocolumn]{SICE21}

\usepackage{amsmath} 
\usepackage{amssymb}  
\usepackage{graphicx,cite}
\usepackage{booktabs, color}
\usepackage{dblfloatfix, float} 
\usepackage{subfigure}

\newtheorem{prob}{Problem}

\title{Temporal Deep Unfolding for Nonlinear Maximum Hands-off Control}

\author{Masako Kishida${}^{1\dagger}$ and Masaki Ogura${}^{2}$}
\speaker{Masako Kishida}

\affils{${}^{1}$Principles of Informatics Research Division, National Institute of Informatics, Tokyo, Japan\\
(E-mail: kishida@nii.ac.jp)\\
${}^{2}$Graduate School of Information Science and Technology, Osaka University, Osaka, Japan\\
(E-mail: m-ogura@ist.osaka-u.ac.jp)\\
}
\abstract{%
This paper proposes a computational technique based on ``deep unfolding'' to solving the finite-time maximum hands-off control problem for discrete-time nonlinear stochastic systems. 
In particular, we seek a sparse control input sequence that stabilizes the system such that the expected value of the square of the final states is small by training a deep neural network. The proposed technique is demonstrated by a numerical experiment.
}

\keywords{%
maximum hands-off control, temporal deep unfolding
}

\begin{document}

\maketitle


\section{Introduction}

Many recent control systems require sparse control inputs. Such control systems include networked control, hybrid vehicles, and railway vehicles. Motivated by these applications, the paradigm of the maximum hands-off control has been introduced in \cite{NagQN16}.

Originally introduced for continuous-time systems \cite{NagQN16}, the maximum hands-off control has been expanded for discrete-time systems \cite{NagOQ16,Rao18,KisNC19}, uncertain systems \cite{KisBN18} and stochastic systems \cite{ExaTT18,ItoIK21}. 
However, there exist few works that deal with nonlinear systems. A possible reason is maybe its computational difficulty, i.e., the optimization problem involves many optimization variables, local minimums, and non-differentiable points, thus it is difficult to solve either analytically or numerically. In this paper, we provide an efficient computational technique to solving a finite-time maximum hands-off control problem for discrete-time nonlinear stochastic systems. 

The idea we employ is so-called \emph{deep unfolding}: which constructs a layer-wise structure by unfolding an iterative algorithm and tunes parameters in the structure such as  step-size and regularization coefficients by standard deep learning techniques. Originally developed to combining the advantages of model-based methods and deep neural networks (DNNs) \cite{HerRW14}, it has been used in signal and image processing and communication systems  \cite{BalS19} as well as average consensus problems \cite{KisOY20,KobOK21}.
To apply this deep unfolding to our control problem, we consider the discrete-time state transition as an iterative algorithm and consider the control inputs as tuning parameters. We call this \emph{temporal deep unfolding}.


\section{Preliminaries}
\subsection{A Short Overview of Maximum Hands-off Control}
The maximum hands-off control aims at minimizing the length of time during which the control input value is nonzero, while achieving given control objectives \cite{NagQN16,NagOQ16}. For a finite-time discrete-time control problem, typically it can be formulated as 
\begin{align}\begin{aligned}
\min_{u} \ \|u\|_0 \\
\text{subject to } & \begin{cases}
\text{system equation}\\
 \text{performance constraint}
 \end{cases}\label{eq:max_ho}
\end{aligned}\end{align}
where $u$ is the vector of all control inputs over the time from 0 to the final time of control input $T-1$ and $\|u\|_0$ denotes  the number of nonzero elements of $u$.
For computational reasons, we usually relax $\|u\|_0$ by $\|u\|_1$, which is the $\ell^1$ norm of $u$ defined by
$\|u\|_1 = \sum_{i=1}^m |u_i|$ where $u_i$ is the $i$-th element of vector $u$ and $m$ is the length of the vector $u$. 
In this paper, we relax $\|u\|_0$ by $\|u\|_p= \left(\sum_{i=1}^m |u_i|^p\right)^{1/p}$ with $p \in ( 0,  1]$ to further seek a sparsity. 
(See \cite{Nag20} for more about the maximum hands-off control.)

\subsection{A Short Overview of Deep Unfolding}
According to \cite{HerRW14}, the idea of deep unfolding can be summarized as ``... given a model-based approach
that requires an iterative inference method, we unfold the iterations into a layer-wise structure analogous to a DNN. We then untie the model parameters across layers to obtain novel neural-network-like architectures that can easily be trained discriminatively using gradient-based methods.”

In short, the deep unfolding can be applied to iterative methods in the form of 
\begin{align}
x_{k+1} = f(x_k, \theta_{k}), \ k = 1, 2,\cdots, K-1, \label{eq:du}
\end{align}
where $\theta_{k}$ is the learning parameter. 
In temporal deep unfolding, we consider the intermediate variables $x_{k}$ as the nodes of layers 1 to $K$ and regarding the equation \eqref{eq:du} as the combination of the transformation and activation function between layers. Then, we can simply train the network with an appropriate loss function and training data to obtain a desired set of parameters $\theta_{k}$.

\section{Finite-time Stabilization by Temporal Deep Unfolding}
Now we are ready to consider the problem of finite-time stabilization and its computational technique. 
\subsection{Problem Setup}
Consider the discrete-time nonlinear stochastic system
\begin{align}
x_{t+1} = f(x_t, u_t, w_t), \ t = 0, 1, 2,\cdots, T-1, \label{eq:sys}
\end{align}
where $x_t \in {\mathbb{R}}^{n}$  is the system state, $u_t \in {\mathbb{R}}^{n_u}$ is the control input, and $w_t \in {\mathbb{R}}^{n_w}$ is the process noise or disturbance, respectively, at discrete time instant $t$.
It is assumed that an initial state $x_0$ is given and the probability distribution of the random variable $w_t$ is known. 

\begin{prob}
For system \eqref{eq:sys}, find a sparse control input sequence $u_0, u_1,\cdots, u_{T-1}$ such that brings $x_T$ near the origin, i.e., 
\begin{align}\begin{aligned}
\min_{u}\ & {{\mathbb{E}}}[\|x_T\|_2^2] +\lambda \|u\|_p\\
\text{subject to }&\eqref{eq:sys}, \label{eq:prob}
\end{aligned}\end{align}
where $\lambda> 0$ and $p \in ( 0,  1]$ are given parameters and $u = [u_0^\top, u_1^\top,\cdots, u_{T-1}^\top]^\top$ is the vector of optimization variables.
\end{prob}

\subsection{Temporal Deep Unfolding for Finite-time Stabilization}

\begin{figure}[h]
    \centering
    \includegraphics[viewport =80 200 670 530, clip, width=\linewidth]{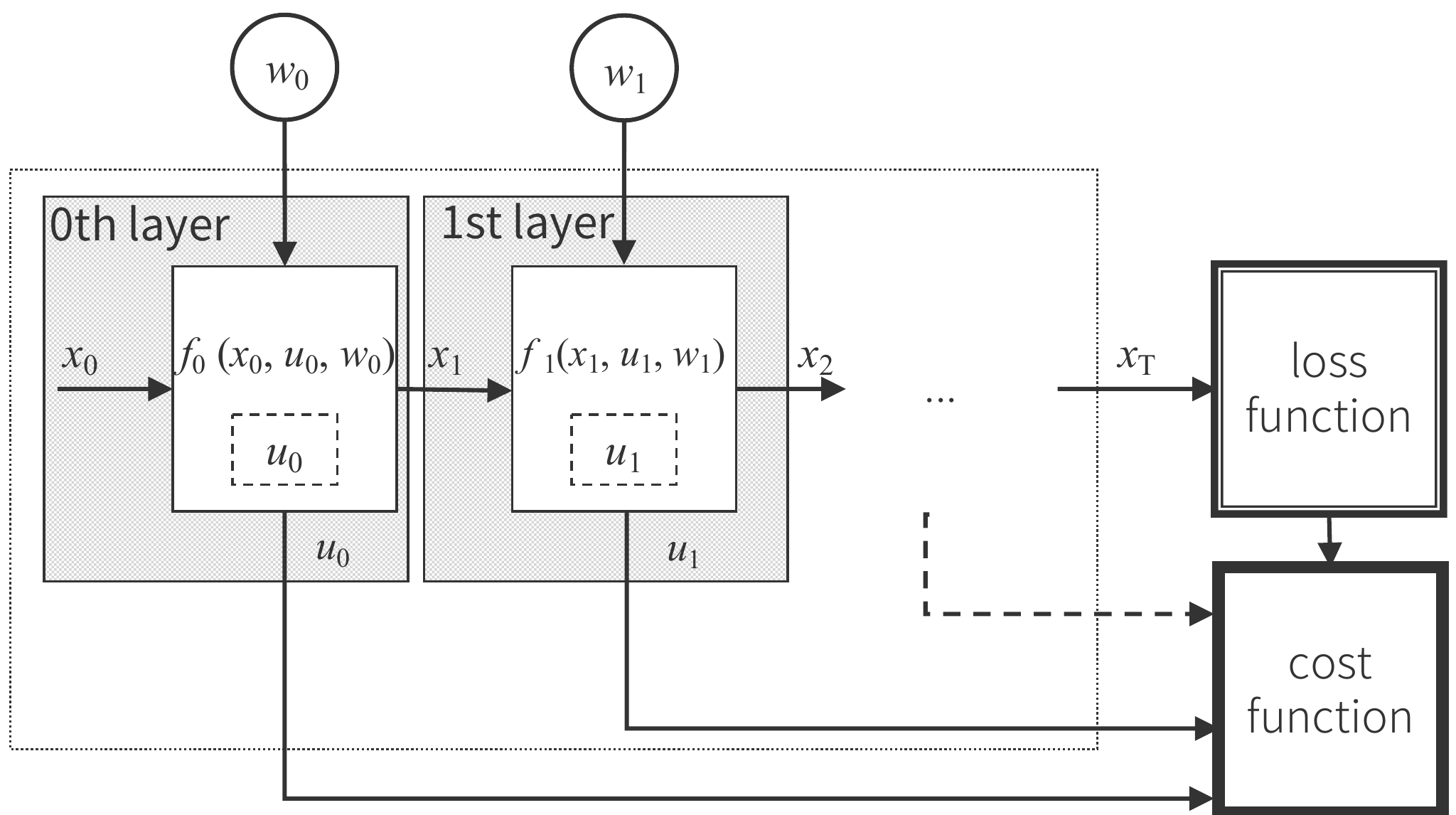}
     \caption{Computation graph structure of the unfolded system, where the learnable parameters are $u_t$ from $t=0$ to $T-1$.} \label{fig:deep_uf} 
  \end{figure}

The proposed architecture of the DNN is shown in Figure \ref{fig:deep_uf}. 
This is obtained by revising the original deep unfolding, i.e., we introduce an input $w_t$ to each layer. 
By this, we can take into account the effect of the disturbance to obtain robust control inputs.

The obtained DNN has $T$ hidden layers with trainable parameters $u = [u_0^\top, u_1^\top,\cdots, u_{T-1}^\top]^\top$.
The output layer computes the loss function: 
\begin{align}
{\mathbb{E}}[\|x_T\|_2^2]. \label{eq:loss}
\end{align}
We also append an $\ell_p$ regularization term:
\begin{align}
\lambda \|u\|_p,
\end{align}
with a hyper-parameter $\lambda>0$, which is given in the problem \eqref{eq:prob}.
Then, the training of the network minimizes the cost function
\begin{align}
{\mathbb{E}}[\|x_T\|_2^2]+ \lambda \|u\|_p.\label{eq:obj}
\end{align}

The training data is $\{w^1, w^2, \cdots, w^I \}$,  where $I$ is the number of training data set
and $w^i = [w_0^{i \top}, w_1^{i \top},\cdots, w_{T-1}^{i \top}]^\top$ is taken from the probability distribution given in the problem.
For any training data, the parameters $u$ are optimized such that  \eqref{eq:obj} is small. Thus unlike standard DNN, the training data is not a set of input-output data, but we have only the input data.

\subsection{Training Techniques}
Here, we briefly explain two techniques that we used during the training to effectively adjust the parameter $u$ that minimizes \eqref{eq:obj} in the numerical experiment in Section \ref{sec:experiments}.

\subsubsection{Incremental Training with Decaying $p$}
In short, the technique of incremental training  \cite{ItoTW19} partitions the original DNN into each layers and trains the partitioned network by gradually increasing the number of layers. In our maximum hands-off control problem, the smaller $p$ tends to yield the sparser control input sequence. 
However, we observed the use of small $p$ from the beginning does not yield an optimal solution. 
Thus, we start with a large $p$ (e.g., $p_0 = 1$) and decrease the value of $p$ gradually as the incremental training proceeds.

More specifically, the first step is to train the $0$th layer. Namely,  the cost function ${\mathbb{E}}\left[\|x_1\|_2^2]+ \lambda \|u_0\|_{p_0}\right]$ is minimized using a number of randomly generated disturbance $w_0$ as the training data $\{w_0^1, w_0^2, \cdots, w_0^{N_p}\}$, where $w_0^i \in {\mathbb{R}}^{n_w}$ are randomly generated from the probability distribution given in the problem. After training the control input $u_0$, we proceed to train the first two control inputs $u_0$ and $u_1$ by appending the first layer to the DNN and decreasing the value of $p$ to $p_1 =  \alpha p_0$. So, the cost function ${\mathbb{E}}\left[\|x_2\|_2^2\right]+ \lambda( \|u_0\|_{p_1} +  \|u_1\|_{p_1})$ is minimized. Here, we use the result from the 0th generation as the initial values for the 0th layer and train the two layers of the DNN. We repeat this process to optimize the all of the control inputs $u_0, \cdots, u_{T-1}$.

\subsubsection{Polishing}
We call a technique that repeats one-shot training while reducing the learning rate, \emph{polishing}. 
We used Adam for stochastic gradient descent. 
Thus, more specifically, by polishing, we mean that the process of repeating the training of all $T$ layers together with the learning rate $(lr)_{i+1} = \beta (lr)_i$ at the $i+1$th repetition. 
We performed polishing after incremental training.

\section{Numerical Experiments} \label{sec:experiments}
In this section, we consider an example of inverted pendulum \cite{Kha02}.
A discrete-time nonlinear model is given by 
\begin{align}
\begin{cases}
x_{t+1} =x_t +\delta y_{t}, \\
y_{t+1}= y_t - \delta \left(\frac{g}{l} \sin x_t +\frac{k}{m} y_t \right) +u_t + w_{t},\\
\end{cases}\label{eq:pendulum}
\end{align}
where $x_t$ is the angle subtended by the rod and the vertical axis through the pivot point ($x_t=0$ is at the bottom equilibrium point and $x_t=\pi$ is the upright equilibrium point),  
$y_t$ is the time derivative of $x_t$, $u_t$ is the control input and $w_t$ is the disturbance. 
We consider swing-up control of the bob from the bottom equilibrium to the upright position.
The meaning of each parameter and its value, as well as other simulation parameters are summarized in Table \ref{table}. 

We first performed incremental training by setting $\lambda = 1$ temporarily and each set of layers is trained using $N_p=2$ noise trajectories. The value of $p$ is decreased as the number of layers increases, starting at the temporal value $p_1=1$ and then $p_{i}= \alpha^{i-1}$ for the training of the first $i$ sets of layers.
After that, we performed polishing using a larger desired value of $\lambda =3\times 10^7$ for 10 times, starting at $(lr)_0 = 1$ and ending at $(lr)_{10} = 0.5^{10}$, using different noise trajectories.
During the polishing, we observed that the control input trajectories smooth out (i.e., small control inputs vanish), while the state trajectories remain nearly the same.

\begin{figure}[ht]
\centering
\subfigure[Control input sequence]
{ \includegraphics[viewport =0 4 400 258, clip, width=.88\linewidth]{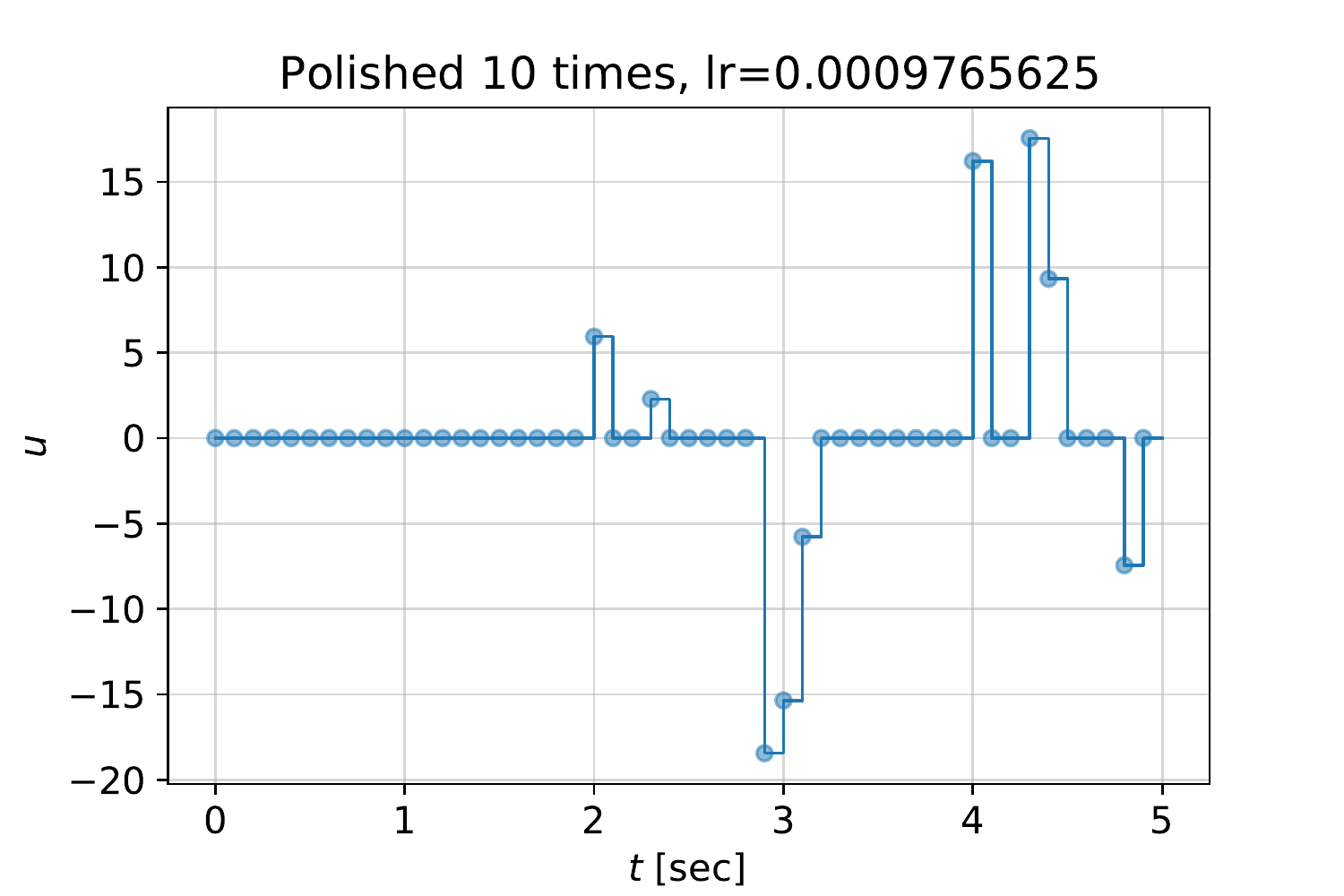}}
\subfigure[State trajectories (the initial state is indicated by a green circle and the target state is indicated by a red circle)]
{ \includegraphics[viewport =0 4 400 258, clip, width=.88\linewidth]{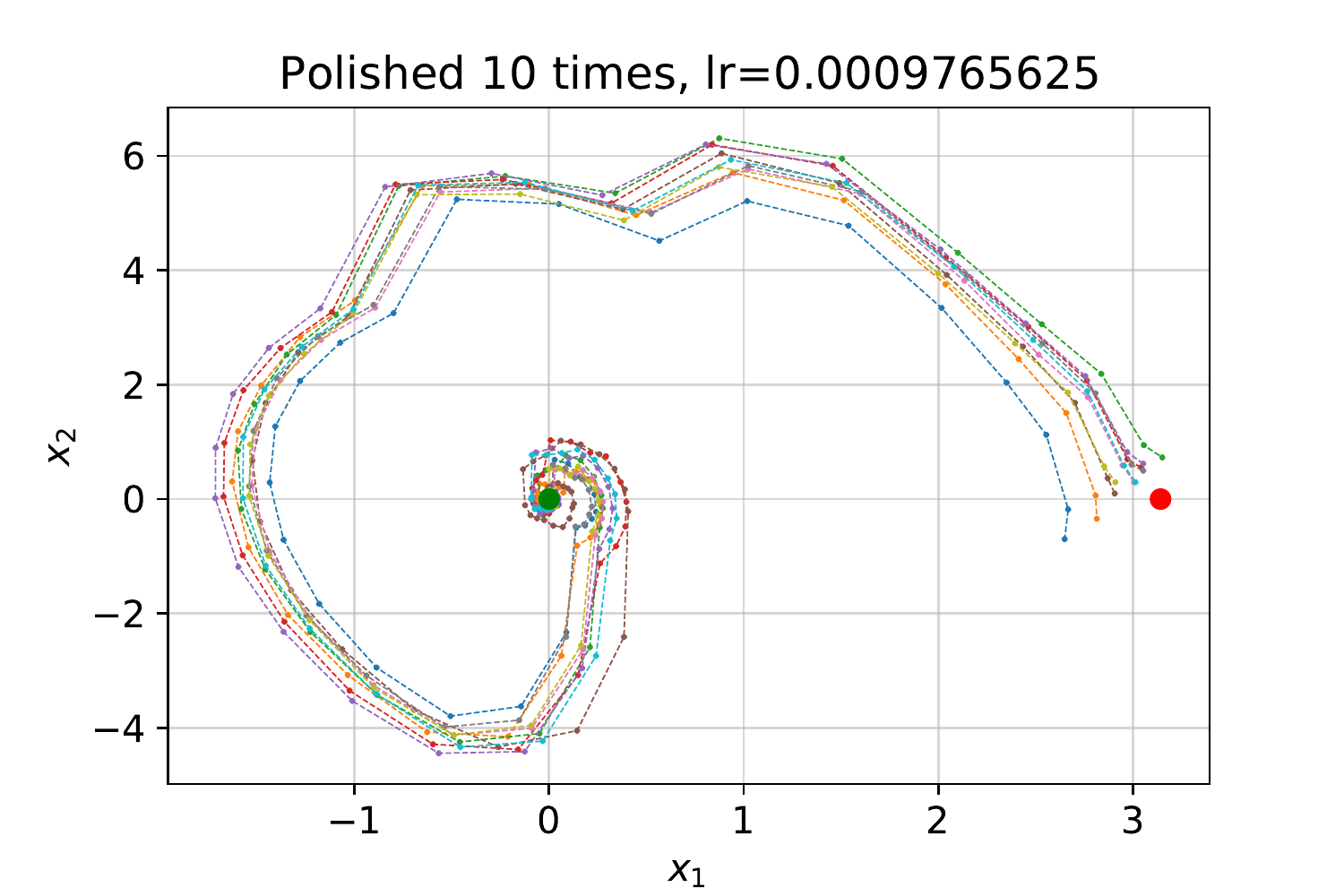}}
\caption{Pendulum system}\label{fig:experiments}\vspace{-.2in}
\end{figure}

The obtained control input sequence is shown in Figure \ref{fig:experiments}(a). 
This control input sequence is then applied to the system \eqref{eq:pendulum} with a randomly generated disturbance sequence over the uniform distribution over $[-1, 1]$. As seen in Figure \ref{fig:experiments}(b), the obtained control input is sparse
and the corresponding state trajectories move toward the target state. 
Note that the control inputs are applied in an open-loop manner, thus we should not expect that the state goes to the target state exactly.

\section{Conclusion} \label{sec:con}
In this paper, we presented how temporal deep unfolding can be applied to nonlinear maximum hands-off control. We also discussed some techniques that are useful when training the DNN to obtaining a desired control input sequence and control performance.

\begin{table}[ht]
\centering
\caption{Parameter values used for the simulation}
\begin{tabular}[t]{lr}
\hline
Parameters&Values\\
\hline
Length of the rod, $l$& 1 \\
Mass of bob, $m$&   1  \\
Friction constant, $k$&  1  \\
Acceleration of gravity, $g$& 9.80665    \\
Sampling time $\delta$ & 0.1\\
Disturbance $w_t$ & uniform i.i.d. $[-1, 1]$\\
Time horizon $T$ & 50\\
Initial state $[x_0,y_0]^\top$ &$ [0, 0]$\\
Final target state $[x_T,y_T]^\top$ &$ [\pi, 0]$\\
Sparsity weight $\lambda$ & $3\times 10^7$\\
Norm $p$ &$1.57\times 10^{-9}$ \\
Number of training data set $N_p$ & 2\\
\quad per layer in incremental training  & \\
Decay rate of $p$, $\alpha$ & 0.667\\
Number of polishing & 10\\
Decay rate of learning rate lr, $\beta$ & 0.5\\
\hline
\end{tabular} \label{table}
\end{table}%

\section*{Acknowledgement}
This work was supported by JST, CREST Grant Number JPMJCR2012, Japan.
\bibliographystyle{IEEEtran}
\bibliography{IEEEabrv,myref}
\end{document}